\documentclass[english]{amsart}

\usepackage{amsmath}
\usepackage{amssymb}
\usepackage{amsthm}
\usepackage{amscd}
\usepackage{babel}
\usepackage[latin1]{inputenc}
\usepackage{graphicx}

\newtheorem{teor}{Theorem}
\newtheorem{cor}{Corollary}
\newtheorem{prop}{Proposition}
\newtheorem{lem}{Lemma}
\theoremstyle{definition}

\newtheorem*{rem}{Remark}
\newtheorem*{conj}{Conjecture}
\newtheorem*{alg}{Algorithm}

\author{José María Grau}
\author{Antonio M. Oller-Marc\'{e}n}
\title{An $\tilde{O}(\log^2(N))$ time  primality test for Generalized Cullen Numbers}
\begin{document}
\maketitle

\begin{abstract}
Generalized Cullen Numbers are positive integers of the form $C_b(n):=nb^n+1$. In this work we generalize some known divisibility properties
of Cullen Numbers and present two primality tests for this family of
integers. The first test is based in the following property of primes from this family: $n^{b^{n}}\equiv (-1)^{b}$ (mod $nb^n+1$). It is
stronger and has less computational cost than Fermat's test (for bases $b$
and $n$) and than Miller-Rabin's test (for base $n$). Pseudoprimes for
this new test seem to be very scarce, only 4 pseudoprimes have been found
among the many millions of Generalized Cullen Numbers tested. We also present a second, more demanding, test for wich no pseudoprimes have been found. This test leads to a ``quasi-deterministic'' test, running in $\tilde{O}(\log^2(N))$ time, which might be very useful in the search of Generalized Cullen Primes.
\end{abstract}

\section{Introduction}
The first major breakthrough in the general theory of primality testing was achieved by Adleman, Pomerance and Rumely (see \cite{ADL}) giving a deterministic primality test running in $\log^{O(\log\log\log n)}$ time. This algorith, later improved by Cohen and Lenstra (see \cite{COH}), is known as the APRCL algorithm.

In 2004 three scholars from Kanpur University (Agrawal, Kayal and Saxena) introduced the AKS algorithm (see \cite{AGR}), which was the first deterministic primality test running in polynomial time. In the second verion of their paper they proved that the running time of their algorithm was $\tilde{O}((\log n)^{7.5})$.  Nevertheless, and despite being one of the cornerstones of Computational Number Theory, this algorithm has not been very useful in practice. This is because numbers for which AKS algorithm is faster than the usual ones are beyond current computation capacity. Even the so-called practical versions of the AKS algorithm (see \cite{BER}, for instance) are not fast enough. As a consequence, prime ``hunters'' focus in families of integers for which primality can be determined by useful algorithms.

For restricted families of integers much faster algorithms are known. The Lucas-Lehmer algorithm (see \cite{LUK}), used for Mersenne Numbers, is deterministic and runs in $\tilde{O}((\log n)^2)$ time. Proth, in \cite{PRO}, gives an algorithm running also in $\tilde{O}(\log n)^2)$ time, which applies to numbers such that $\nu_2(n-1)>\frac12\log_2 n$ where $2^{\nu_2(m)}$ is the biggest power of 2 dividing $m$ and provided an integer $a$ is given such that the Jacobi symbol $\left(\frac{a}{n}\right)=-1$. Proth's algorithm is not deterministic for every $n$. Later, Williams \cite{WIL} or Konyagin and Pomerance \cite{KON} have extended these techniques to wider families of integers. 

Positive integers of the form $n2^{n}+1$ are called Cullen Numbers and were first introduced by Father James Cullen in 1905 (see \cite{CUL}, \cite[B20]{GUY} or \cite{LUC} for instance). Primes of this form are very scarce (in fact, in \cite{HOO} it is shown that almost all Cullen Numbers are composite). Primality criteria suitable for Cullen Numbers have been presented and discussed in \cite{ROB}. The only known Cullen Primes are those for $n$ equal to:
\begin{center}
1,\ 141,\ 4713,\ 5795,\ 6611,\ 18496,\ 32292,\ 32469,\ 59656,\ 90825,\ 262419, \\
361275,\ 481899,\ 1354828,\ 6328548,\ 6679881\text{ (sequence A005849 in OEIS)}.
\end{center}
The largest known Cullen prime is $6679881\times 2^{6679881}+1$. It is a megaprime with 2,010,852 digits and was discovered by a \textit{PrimeGrid} participant from Japan. It is the fifteenth biggest known prime. In \cite{BYF} a ``quasi-deterministic'' test for Cullen Numbers was given. 

A quite straightforward generalization of these numbers are the so-called Generalized
Cullen Numbers (GCN for short) which are integers of the form $C_{b}(n):=nb^{n}+1$. This family was introduced by H. Dubner
in \cite{DUB} and is one of the main sources for prime number ``hunters''. There exists a distributed computing project (http://www.primzahlenarchiv.de/) to find Generalized Cullen
Primes (GCP for short) with the biggest GCP being $C_{151}(139948)$, an integer with 304,949 digits. Noteworthy, for 29 values of $b$
smaller than 200 no GCP has been found.

To date, no specific primality test for Generalized Cullen Numbers has been introduced. This is the main goal of this work. The paper is organized as follows. In the second section we generalize some known divisibility properties of Cullen Numbers. In Section 3 we present two probabilistic primality tests for GCN. The first test (TEST1) is based in the fact that $n^{b^{n}}\equiv (-1)^{b}$ (mod $C_{b}(n)$) for every GCP $C_b(n)$. This
test is stronger and has less computational cost than Fermat's test (for
bases $b$ and $n$) and than Miller-Rabin's test (for base $n$) and seems to have very few pseudoprimes. Thus, the probability of error is extremely small: among the millions of numbers tested, only 4 pseudoprimes have been found. We also present another test (TEST2), more demanding than TEST1, for which no pseudoprime has been found. In the fourth section we present a ``quasi-deterministic'' version of TEST2, which has allowed to certify the primality of nearly every known GCP with the use of very modest technological resources in just a few minutes. We are convinced that, with the use of better technology directed to an efficient modular exponentiation, this algorith would help to break records for GCP. Finally, in Section 5, we stablish the computational complexity of the presented tests, we give the running time of our algorithm for various cases and we close the paper with an important conjecture.

\section{Some divisibility properties}
Although the main goal of the paper is to present primality tests for Generalized Cullen Numbers, it is interesting to study some divisibility properties of such numbers. First of all, we are interested in finding families of composite Generalized Cullen Numbers. A first result in this direction goes as follows.

\begin{prop}
Let $n_b(k,p)=(b^k-k)(p-1)-k$ and let $p$ be a prime not dividing $b$. Then $p$ divides $C_b(n_b(k,p))$.
\end{prop}
\begin{proof}
It is clear that $n_b(k,p)\equiv -b^k\ \textrm{(mod $p$)}$. Now, since $b^{p-1}\equiv 1\ \textrm{(mod $p$)}$ we have that:
$$C_b(n_b(k,p))\equiv -b^{(b^k-k)(p-1)}+1\equiv 0\ \textrm{(mod $p$)}.$$
\end{proof}

Observe that, if $p$ divides $C_b(n)$ then $p$ does not divide $b$, so if $n\neq n_b(k,p)$, we can apply the previous proposition to find another composite Generalized Cullen Number with the same base. Nevertheless, this process can be applied only in one step; i.e., given a prime divisor of $C_b(n)$ we only find (at most) another $m$ such that $p$ also divides $C_b(m)$. Thus, it is interesting to find a process that allows us to construct infinite families of Generalized Cullen Numbers divisible by the same prime. The next proposition goes in this direction.

\begin{prop}
Let $p$ be a prime dividing $C_b(n)$ and let $h_p=exp_p(b)$; i.e., the smallest integer such that $b^{h_p}\equiv 1\ \textrm{(mod $p$)}$. Then $p$ also divides $C_b(n+mph_p)$ for every integer $m$.
\end{prop}
\begin{proof}
If $m=1$ we have:
$$C_b(n+ph_p)=b^{ph_p}(C_b(n)-1)+ph_pb^{n+ph_p}+1\equiv nb^n+1=C_b(n)\equiv 0\ \textrm{(mod $p$)}$$
and the result follows inductively.
\end{proof}

The propositions above generalize known results for the case $b=2$ that can be found in \cite{KEL} and \cite{CUN}.

Now, given two Generalized Cullen Numbers, it can be interesting to study their common divisors. For example, if we consider $C_b(n)$ and $C_{\beta}(n)$ it is easy to see that any common divisor of these numbers must also divide $|b^n-\beta^n|$. If we restrict ourselves to Generalized Cullen Numbers of the same base we can present a more interesting result.

\begin{prop}
Let $C_b(n)$ and $C_b(m)$ be two different Generalized Cullen Numbers with $\alpha n=\beta m$. If $d$ is a common divisor of $C_b(n)$ and $C_b(m)$, then $d$ also divides $|n^{\alpha}+(-1)^{\alpha+\beta-1}m^{\beta}|$.
\end{prop}
\begin{proof}
First of all note that $d$ cannot divide $b$. Also, since $d$ is a common divisor, it must, assuming $n\geq m$, divide $|C_b(n)-C_b(m)|=b^n|n-mb^{m-n}|$. Consequently $d$ divides $|n-mb^{m-n}|$.

Now, $d$ also divides $nC_b(n)=n^2b^n+n$ and it follows that $d$ divides $|n^2b^n+mb^{m-n}|$. If, for instance, $n<m-n$ we get that $d$ divides $|n^2-mb^{m-2n}|$. If $m-n<n$ we would get that $d$ divides $|m+n^2b^{2n-m}|$.

Clearly we can proceed in this way until both powers of $b$ are the same. Furthermore, we will need to perform the previous computations exactly $\alpha+\beta-1$ times and at every step the middle sign changes and the exponent of either $m$ or $n$ increases by 1. Thus, the desired result is finally obtained. We omit the details.
\end{proof}

It is worth remarking that the previous proposition does not depend on $b$. Thus, if $A_{m,n}=|n^{\alpha}+(-1)^{\alpha+\beta-1}m^{\beta}|$ is a prime, then $C_b(n)$ and $C_b(m)$ are either coprime or their greatest common divisor is $A_{m,n}$ for every value of $b$.

\section{Two probabilistic primality tests}
This section is devoted to present two probabilistic primality tests for Generalizados Cullen Numbers. The first one will be compared with Fermat and Miller-Rabin for some witnesses. The second one will be the basis of a ``quasi-deterministic'' test which will we introduced in the next section.

The proposition below presents the property of GCP in which TEST1 will be based.

\begin{prop}
If $C_b(n)$ is prime, then $n^{b^n}\equiv (-1)^b$ (mod $C_b(n)$).
\end{prop}
\begin{proof}
Clearly, $nb^n\equiv -1$ (mod $C_b(n)$). On the other hand, since $b$ and $C_b(n)$ are coprime, we have that $b^{nb^n}\equiv 1$ (mod $C_b(n)$).

Now, taking this into account:
$$-(-1)^b\equiv (-1)^{b^n-1}\equiv (nb^n)^{b^n-1}\equiv n^{b^n-1}b^{nb^n-n}\equiv n^{b^n-1}b^{-n}\ \textrm{(mod $C_b(n)$)},$$
from where it follows that:
$$n^{b^n}\equiv -(-1)^bnb^n\equiv (-1)^b\ \textrm{(mod $C_b(n)$)},$$
where negative exponents make sense since we are working over a field.
\end{proof}

Let us relate now this test with Fermat and Miller-Rabin primality tests. In fact, we will se that our test is stronger than both of them in the sense that if a Generalized Cullen Number passes our test, it will also pass Fermat and Miller-Rabin tests for certain choices of the base.

\begin{prop}
If $n^{b^n}\equiv (-1)^b$ (mod $C_b(n)$) for an odd $C_b(n)$, then $C_b(n)$ is a Fermat (or weak) probable prime to base $n$.
\end{prop}
\begin{proof}
We have that:
$$n^{C_b(n)-1}=n^{nb^n}=(n^{b^n})^n\equiv (-1)^{bn}=1\ \textrm{(mod $C_b(n)$)},$$
since $bn$ must be even.
\end{proof}

\begin{rem} 
Although TEST1 is very similar to Fermat's test to base
$n$, it has turned to be more subtle. We have only found 4 pseudoprimes for our test. Namely: $C_{80}(2)=12801$, $C_{3570}(3)=136497879001$, $C_{570}(4)=422240040001$ and $C_{1470}(4)=18677955240001$). On the other hand, three more pseudoprimes ($C_{7}(4)$, $C_{63336}(2)$ and $C_{2355990}(2)$) appear for Fermat's test to base $n$. Observe that $C_{1470}(4)$ and $C_{570}(4))$ are Carmichael Numbers.
\end{rem}

\begin{prop}
If $n^{b^n}\equiv (-1)^b$ (mod $C_b(n)$), then $C_b(n)$ is a Fermat (or weak) probable prime to base $b$.
\end{prop}
\begin{proof}
Since $b$ and $n$ are coprime with $C_b(n)$, they both have an inverse modulo $C_b(n)$. Moreover, $n^{-1}\equiv -b^n$ (mod $C_b(n)$).

Now:
$$(-1)^b\equiv n^{b^n}\equiv (-b^{-n})^{b^n}\equiv (-1)^bb^{-nb^n}\ \textrm{(mod $C_b(n)$)}.$$
Thus, $b^{-nb^n}\equiv 1$ (mod $C_b(n)$) and, consequently, $b^{C_b(n)-1}=b^{nb^n}\equiv 1$ (mod $C_b(n)$).
\end{proof}

\begin{rem} 
Although TEST1 is theoretically stronger than Fermat's test to base $b$, we have not found
pseudoprimes for this test which are not pseudoprimes for TEST1.
\end{rem}

\begin{prop}
If $n^{b^n}\equiv (-1)^b$(mod $C_b(n)$), then $C_b(n)$ is a strong probable prime to base $n$; i.e., it passes Miller-Rabin primality test.
\end{prop}
\begin{proof}
Put $n=2^rh$ and $b=2^sk$ with $k,h$ odd integers and $r,s\geq0$. Then, $C_b(n)=2^rh(2^sk)^{n}+1=2^{t}m+1$ with $t=r+ns$ and $m=hk^{n}$. We have that $n^{b^n}\equiv (-1)^b$ (mod $C_b(n)$). Consequently $n^{2^tm}=n^{nb^n}\equiv (-1)^{nb}\equiv 1$ (mod $C_b(n)$), since $nb$ must be even. If $n^m\not\equiv 1$ (mod $C_b(n)$) it follows that $n^{2^{t-1}m}\equiv -1$ (mod $C_b(n)$) and $C_b(n)$ passes Miller-Rabin test to base $n$.
\end{proof}

The following result is stronger than Proposition 4 and will give rise to another probabilistic test that we will denote by TEST2. This test involves cyclotomic polynomials of prime index and is more demanding than TEST1. It will be the basis of a ``quasi-deterministic'' test in the next section.

\begin{teor} 
Let $p$ be a prime number and $b=p^{m}b'$ with $p$ not dividing $b'$. If $C_{b}(n):=nb^{n}+1$ is prime, then one of the following holds:
\begin{itemize}
\item[i)] $(-n)^{\frac{b^{n}}{p^{i}}}\equiv1$ (mod $C_{b}(n)$) for every $i\in\{0,...,nm\}$.

\item[ii)] There exists $K<mn$ such that $(-n)^{\frac{b^{n}}{p^{i}}}\equiv 1$ (mod $C_{b}(n)$) for every $i\in\{0,...,K\}$ and $\Phi_{p}((-n)^{\frac{b^{n}}{p^{K+1}}})\equiv 0$ (mod $C_b(n)$), where $\Phi_p$ is the $p$-th cyclotomic polynomial.
\end{itemize}
\end{teor}
\begin{proof} 
Proposition 4 implies that $(-n)^{\frac{b^{n}}{p^{0}}}\equiv1$ (mod $C_b(n)$). If i) does not hold, let $K<mn$ be the biggest integer such that $(-n)^{\frac{b^{n}}{p^{K}}}\equiv1$ (mod $C_b(n)$). Put $x=(-n)^{\frac{b^{n}}{p^{K+1}}}$. Then $0\equiv x^p-1=(x-1)\Phi_p(x)$ (mod $C_b(n)$). The maximality of $K$ implies that $x-1\not\equiv 0$ (mod $C_b(n)$) so, since $C_b(n)$ is prime, $\Phi_p(x)\equiv 0$ (mod $C_b(n)$) and the proof is complete.
\end{proof}

Every Generalized Cullen Number satisfying Theorem 1 for some $p$, prime divisor of $b$, will be certified as a probable prime for TEST2 to base $p$. If it is composite we will name it as a pseudoprime for TEST2 to base $p$. Among the four pseudoprimes found for TEST1, the only one that passes TEST2 for some prime divisor of $b$ is the Carmichael number $C_{1470}(4)$ which is a pseudoprime for TEST2 to base 2. Nevertheless, $C_{1470}(4)$ is certified as composite since it does not pass TEST2 for any other prime divisor of 1470. The other three pseudoprimes for TEST1 do not pass TEST2 for any base. The authors have not found any composite Generalized Cullen Number passing TEST2 for every prime divisor of $b$. In fact they conjecture that such GCN does not exist.

\section{A ``quasi-deterministic'' test}

We will now see that passing TEST2, together with a bounding condition on $K$, gives a sufficient condition for primality. 

\begin{teor}
Let $b=p_1^{r_1}\cdots p_s^{r_s}$ be a positive integer. If there exist $p_j$ prime dividing $b$ and $K\leq nr_j$ such that:
\begin{itemize}
\item [i)] $\Phi_{p_j}((-n)^{\frac{b^n}{p_j^{K+1}}})\equiv 0$ (mod $C_b(n)$).
\item[ii)] $nr_j-K>\frac12\log_{p_j} nb^n=\frac12\log_{p_j}n+\frac{n}{2}\log_{p_j}b$.
\end{itemize}
Then $C_b(n)$ is prime.
\end{teor}
\begin{proof}
If $q$ is a prime divisor of $C_b(n)$ we have that $\Phi_{p_j}((-n)^{\frac{b^n}{p_j^{K+1}}})\equiv 0$ (mod $q$). It follows that the order of $(-n)^{\frac{b^n}{p_j^{K+1}}}$ in $\mathbb{Z}^{*}_{q}$ is exactly $p_j$ and, consequently, the order of $-n$ is a divisor of $\frac{b^n}{p_j^K}=p_1^{nr_1}\cdots p_j^{nr_j-K}\cdots p_s^{nr_s}$. Moreover, $(-n)^{\frac{b^n}{p_j^{K+1}}}$ is not congruent with 1. For if it was, then $0\equiv \Phi_{p_j}(1)$ (mod $q$) which is a contradiction being $q$ and $p_j$ coprime. Thus, the order of $-n$ does not divide  $\frac{b^n}{p_j^{K+1}}=p_1^{nr_1}\cdots p_j^{nr_j-K-1}\cdots p_s^{nr_s}$. As a consequence, the order of $-n$ is a multiple of $p_j^{nr_j-K}$ and it follows that $p_j^{nr_j-K}|q-1$. Finally we obtain that:
$$q\geq p_j^{nr_j-K}+1>p_j^{\frac12\log_{p_j} nb^n}+1=\sqrt{nb^n}+1\geq\sqrt{C_b(n)}.$$
This must hold for every $q$ prime divisor of $C_b(n)$. Clearly it is a contradiction unless $q=C_b(n)$ is its only prime divisor and the result follows.
\end{proof}

If $b$ is a prime-power the previous result can be slightly simplyfied in the following way.

\begin{cor}
Let $b=p^m$ with $p$ a prime. If $\Phi_{p}((-n)^{p^{K-1}})\equiv 0$ (mod $C_b(n)$) with $nm\geq K>\frac{mn}{2}+\frac12\log_p(n)$, then $C_b(n)$ is prime.
\end{cor}

\begin{rem}
The previous Corollary generalizes \cite[Teorema 4]{BYF} for $b=2$ and $m=1$. Observe that in such case $\Phi_{2}((-n)^{2^{K-1}})=n^{2^{K-1}}+1$.  
\end{rem}

\section{Computational complexity}
Let us present the pseudo code of an algorithm implementing TEST1 and TEST2. We will also justify its polynomial complexity.

\begin{alg} (TEST1 and TEST2 to base  $P$)

INPUT: $n$ and  $b=p_{1}^{m_{1}}\ldots p_{s}^{m_{s}};$ $M:=m_{j};$ $P:=p_{j}$
and $N:=nb^{n}+1$

\textbf{Step 1}: If $n^{b^{n}}\not\equiv(-1)^{b}$ (mod $N$)

then RETURN: $N$  is a COMPOSITE NUMBER. Stop.

\textbf{Step 2}: If $n^{b^{n}}\equiv(-1)^{b}$ (mod $N$)

then RETURN: $N$ is a PROBABLE PRIME for TEST1

and go to Step 3.

\textbf{Step 3}: Compute $K:= \max\{i\leq nM\ |\ ((-n)^{\frac{b^{n}}{P^{i}}}\equiv1\ (\textrm{mod}\ N)\}$

\textbf{Step 4}: If $K=nM$

then  RETURN: $N$ is a PROBABLE PRIME for TEST2 to base $P$. Stop.

\textbf{Step 5}: If $K<nM$ and $\Phi_{P}((-n)^{\frac{b^{n}}{P^{K+1}}})\equiv 0$ (mod $N$)

then RETURN: $N$ is a PROBABLE PRIME for the TEST2 to base $P$

and go to Step 7.

\textbf{Step 6}: If $K<nM$ and $\Phi_{P}((-n)^{\frac{b^{n}}{P^{K+1}}})\neq 0$ (mod $N$)

then RETURN: $N$ is a COMPOSITE NUMBER. Stop.

\textbf{Step 7}:  If $nM-K>\frac{1}{2}\log_{P}(N-1)$

then RETURN: $N$ is a PRIME NUMBER. Stop.
\end{alg}

The correctness of the algorithm is a straightforward consequence of Proposition 4 and Theorems 1 and 2. To study its complexity, we first need to present a technical lemma. Its proof is elementary and we omit.

\begin{lem}
Let $A=\{a_0, a_1,\dots, a_n\}$ a set with $a_n\neq a_0$ and with the property $a_s\neq a_0\Rightarrow a_{s+1}\neq a_0$. If the computation of each $a_i$ is of complexity $O(h)$, then the complexity of computing $\max\{i\ |\ 0\leq i\leq n,\ a_i=a_0\}$ is $O(h\log_2(n))$.
\end{lem}

\begin{teor}
If $N:=nb^{n}+1$, the complexity of the algorithm above is $$O(\log^{3}(W(N))\log^2(N))\subset \widetilde{O}(\log^{2}(N)),$$ where $W$ is Lambert's W function.
\end{teor}
\begin{proof}
Complexity of steps 1 and 2 is that of the modular exponentiation $n^{(N-1)/n}$ (mod $N$). Taking into account that $O(n)=O(W(N))\subset O(\log(N))$  and that products modulo $N$ can be performed by Schoenhage-Strassen algorithm (see \cite{SCH}) with complexity
$$O(\log(N) \log(\log(N))\log(\log(\log(N)))),$$
we get the complexity of these two steps to be:
$$T(N):=O(\log\frac{N}{W(N)}\log(N) \log(\log(N))\log(\log(\log(N)))).$$

Step 3 requires to perform at most $\log(nM)$ modular exponentiations with complexity analogous to that of step 1. Consequently, by the lemma above, the total complexity will be

 $$O(T(N)\log(W(N))\subset O(\log^{3}(W(N))\log^2(N))\subset \widetilde{O}(\log^{2}(N)).$$

Steps 4 through 7 do not increase the complexity since they are mere verifications of equalities and inequalities.
\end{proof}

From a computational point of view, this work is promising because it presents primality tests for GCN whose computational complexity is
the same as that of modular exponentiation. This is the case, since only a relatively small number of modular power has to be computed. Other tests of general character are clearly inferior. For example, Lucas test, which seems appropriate here due to the easy factorization of $C_{b}(n)-1$ has problems if $n$ has many prime divisors. Altough the techology used by the authors limits to \textit{PowerMod} command in $\textrm{Mathematica}^{\texttt{TM}}$ 6.0 (in an Intel core2 Duo P7450 @ 2.13 GHz with 4Gb of RAM), they have been able to certify primality for nearly every known GCP. We are sure that the tests presented in this work, using a better technology focused on efficient modular exponentiation, will allow to break records in the family of GCN.

Just to enlighten what we have just said, let us compare our primality test with the \textit{PrimeQ} command implemented in $\textrm{Mathematica}^{\texttt{TM}}$. According to $\textrm{Mathematica}^{\texttt{TM}}$ manual, this command uses the multiple Rabin-Miller test in bases 2 and 3 combined with a Lucas pseudoprime test. Below we present the runnig times for the certification of some known GCP (for  $b=3$, $8$ and $20$). We also show the number of modular exponentiations performed in the third step of the algorithm.
$$
\begin{tabular}{|l|l|l|l|l|l|l|l|l|l|l|l|l|l|}
\hline
b=3, n = &  1400 & 1850 & 2848 & 4874 & 7268 &19290 & 337590 \\ \hline
Number of digits of $C_{3}(n)$ & 672&886&1363&2330 & 3472 &  9208 &161077 \\ \hline
$K+1$ (Step 3) &  1 & 2 & 2 & 1 & 1 & 1 & ?
\\ \hline

\textbf{TEST2 to base 3} (Time in s.)  & 0.06 & 0.17 & 0.56 & 1.49 &
4.55 & 54.16 & ?  \\ \hline
\textbf{PrimeQ }(Time in s.)   & 0.09 & 0.17 & 0.50 & 2.04 &
6.06 & 71.19 & ?  \\ \hline
\end{tabular}
$$

$$
\begin{tabular}{|l|l|l|l|l|l|l|l|}
\hline
b=8, n= & 5 & 17 & 23 & 1911 & 20855 & 35945 & 42816 \\ \hline
Number of digits of $C_{8}(n)$ &6&17&23 & 1730 & 18839 & 32467 &38672 \\ \hline
$K+1$ (Step 3) & 2 & 2 & 2 & 3 & 6 & 5 & 2
\\ \hline
 \textbf{TEST2 to base 2} (Time in s.) &  &  &  & 0.87 & 788. & 2811. & 2112. \\ \hline
\textbf{PrimeQ}(Time in s.) &  &  &  & 0.93 & 464. & 1933. & 2930. \\ \hline
\end{tabular}
$$

$$
\begin{tabular}{|l|l|l|l|l|l|l|l|}
\hline
b=20, n= & 3 & 6207 & 8076 & 22356  \\ \hline
Number of digits of $C_{20}(n)$  & 5 & 8080 & 10512 &29091 \\ \hline
 $K+1$ (Step 3) & 1 & 1 & 1 &2 \\ \hline
 \textbf{TEST2 to base 5} (Time in s.) & 0. &32.1& 62.4 & 1347.  \\ \hline
 \textbf{PrimeQ} (Time in s.) & 0. & 50.88 &99.6 & 1396.  \\ \hline
\end{tabular}
$$

Note that \textit{PrimeQ} is a probabilistic test (even though according to $\textrm{Mathematica}^{\texttt{TM}}$ manual no pseudoprime for this test has still been found). In the examples above our test certified primality in a deterministic way and faster (except for two cases) than \textit{PrimeQ}. The deterministic version of \textit{PrimeQ} implemented in$\textrm{Mathematica}^{\texttt{TM}}$, \textit{ProvablePrimeQ}, was so slow that it is not worth presenting the comparison.

To finish, we have to point out that the presented algorithm is not strictly a deterministic primality test. We could say it is a ``quasi-deterministic'' test in the sense that if the output is PRIME NUMBER, the tested integer is certainly a prime, but it could be possible that a GCP was not certified as such. Nevertheless, computational evidence suggest that this is not the case for moderately big values of $n$; in fact, such primes have be found only for $n<b$. Moreover, the experiments also suggest that in step 3 the value of $K$ is always very small with respect to $n$. This leads us to the following conjecture.

\begin{conj}
If $n>b$, then the algorithm is deterministic.
\end{conj}

\section*{Acknowledgememts}
We are grateful to L. M. Pardo Vasallo for his help in computational complexity aspects. We are also grateful to P. Berrizbeitia and J.G. Fernandes for providing the prepint of their paper \cite{BYF} that inspired this work.

\bibliography{./refcullen}
 \bibliographystyle{plain}
\end{document}